\documentclass{proc-l}

\usepackage{epsf}
\usepackage{epsfig}
\usepackage{graphicx}
\usepackage{xcolor}
\usepackage{bm}
\usepackage{tabularx}
\usepackage{amstext}
\usepackage{float}
\usepackage{array, nicefrac, caption}

\newcommand{\divprop}{\lfloor}
\newtheorem{lemma}{Lemma}
\newtheorem{theorem}{Theorem}
\newtheorem{corollary}{Corollary}
\newtheorem{definition}{Definition}

\newtheorem{proposition}{Proposition}

\allowdisplaybreaks
\begin{document}

\frenchspacing

\title{Recursively divisible numbers}
\author{Thomas Fink}
\address{London Institute for Mathematical Sciences, Royal Institution, 21 Albemarle St, London W1S 4BS, UK}

\date{\today}

\begin{abstract}
We introduce and study the recursive divisor function, a recursive analog of the usual divisor function:
$\kappa_x(n) = n^x + \sum_{d\lfloor n} \kappa_x(d)$,
where the sum is over the proper divisors of $n$.
We give a geometrical interpretation of $\kappa_x(n)$, 
which we use to derive a relation between $\kappa_x(n)$ and $\kappa_0(n)$.
For $x \geq 2$, we observe that $\kappa_x(n)/n^x < 1/(2-\zeta(x))$.
We show that, for $n \geq 2$, $\kappa_0(n)$ is twice the number of ordered factorizations, a problem much studied in its own right.
By computing those numbers that are more recursively divisible than all of their predecessors,
we recover many of the numbers prevalent in design and technology, and suggest new ones which have yet to be adopted.
\vspace*{-0.25in}
\end{abstract}

\maketitle

\section{Introduction}
\subsection{Recursive divisor function} \noindent 
In this paper we introduce and study the recursive divisor function:
\begin{definition}	\label{adef}
	\begin{equation*}
		\kappa_x(n) = n^x + \sum_{d\divprop n} \kappa_0(d),
		\label{MainDef}
	\end{equation*}
where $d \divprop n$ means $d \vert n$ and $d < n$.
\end{definition}
\noindent
It can be thought of as the recursive analogue of the usual divisor function:
\begin{equation}
	\sigma_x(n) = \sum_{d|n} d^x.
	\label{AN}
\end{equation}
The recursive divisor function considers not only the divisors of $n$,
but also the divisors of the resultant quotients, and the divisors of those resultant quotients, and so on.
For example, 
$\kappa_0(10) = 1 + \kappa_0(1) + \kappa_0(2) + \kappa_0(5) = 6$, and 
$\kappa_1(10) = 10 + \kappa_1(1) + \kappa_1(2) + \kappa_1(5) = 20$.
Values of $\kappa_0, \kappa_1$ and $\kappa_2$ are given in Table 1.
\subsection{Example} \noindent
Consider one of the earliest references to a number that can be divided into equal parts in many ways.
Plato writes in his \emph{Laws} that the ideal population of a city is 5040, since this number has more divisors than any number less than it.
He observes that 5040 is divisible by 60 numbers, including one to 10.
A highly divisible population is useful for dividing the city into equal-sized sectors for administrative, social and military purposes.
\\ \indent
This conception of divisibility can be extended.
Once the city is divided into equal parts, it is often necessary to divide a part into equal subparts.
For example, 	if 5040 is divided into 15 parts of 336, each part can in turn be divided into subparts in 20 ways, since 336 has 20 divisors.
But			if 5040 is divided into 16 parts of 315, each part can be divided into subparts in only 12 ways, since 315 has 12 divisors.
Thus the division of the whole into 15 parts offers more optionality for further subdivisions than the division into 16 parts.
Similar reasoning can be applied to the division of the subparts into sub-subparts, and so on, in a recursive way.
\begin{table}[b!]
\begin{small}
\begin{tabular*}{\textwidth}{@{\extracolsep{\fill}}rrrrrrrrrrrr}
$n$ & $\kappa_0$ 	&	$\kappa_1$ 	& 	$\kappa_2$ &  \mbox{\hspace{0.26in}}	
$n$ & $\kappa_0$ 	& 	$\kappa_1$ 	& 	$\kappa_2$ &  \mbox{\hspace{0.26in}}	
$n$ & $\kappa_0$ 	& 	$\kappa_1$ 	& 	$\kappa_2$ \\	
 1 & 1 & 1 & 1 & 21 & 6 & 34 & 502 & 41 & 2 & 42 & 1682 \\
 2 & 2 & 3 & 5 & 22 & 6 & 38 & 612 & 42 & 26 & 132 & 2636 \\
 3 & 2 & 4 & 10 & 23 & 2 & 24 & 530 & 43 & 2 & 44 & 1850 \\
 4 & 4 & 8 & 22 & 24 & 40 & 116 & 992 & 44 & 16 & 106 & 2698 \\
 5 & 2 & 6 & 26 & 25 & 4 & 32 & 652 & 45 & 16 & 96 & 2416 \\
 6 & 6 & 14 & 52 & 26 & 6 & 44 & 852 & 46 & 6 & 74 & 2652 \\
 7 & 2 & 8 & 50 & 27 & 8 & 46 & 832 & 47 & 2 & 48 & 2210 \\
 8 & 8 & 20 & 92 & 28 & 16 & 74 & 1114 & 48 & 96 & 304 & 4088 \\
 9 & 4 & 14 & 92 & 29 & 2 & 30 & 842 & 49 & 4 & 58 & 2452 \\
 10 & 6 & 20 & 132 & 30 & 26 & 104 & 1388 & 50 & 16 & 112 & 3316 \\
 11 & 2 & 12 & 122 & 31 & 2 & 32 & 962 & 51 & 6 & 74 & 2902 \\
 12 & 16 & 42 & 234 & 32 & 32 & 112 & 1520 & 52 & 16 & 122 & 3754 \\
 13 & 2 & 14 & 170 & 33 & 6 & 50 & 1222 & 53 & 2 & 54 & 2810 \\
 14 & 6 & 26 & 252 & 34 & 6 & 56 & 1452 & 54 & 40 & 190 & 4392 \\
 15 & 6 & 26 & 262 & 35 & 6 & 50 & 1302 & 55 & 6 & 74 & 3174 \\
 16 & 16 & 48 & 376 & 36 & 52 & 176 & 2196 & 56 & 40 & 196 & 4672 \\
 17 & 2 & 18 & 290 & 37 & 2 & 38 & 1370 & 57 & 6 & 82 & 3622 \\
 18 & 16 & 54 & 484 & 38 & 6 & 62 & 1812 & 58 & 6 & 92 & 4212 \\
 19 & 2 & 20 & 362 & 39 & 6 & 58 & 1702 & 59 & 2 & 60 & 3482 \\
 20 & 16 & 58 & 586 & 40 & 40 & 156 & 2464 & 60 & 88 & 346 & 6318 
\end{tabular*}
\end{small}
\vspace{0.1in}
\caption{\small 
{\bf Values of} 
\boldmath $\kappa_0,$ \unboldmath
\boldmath $\kappa_1$ \unboldmath
{\bf and} 
\boldmath $\kappa_2.$ \unboldmath
A Mathematica algorithm for $\kappa_0(n)$ is: 
{\tt 
n = 1; $\kappa$ = \{\}; 
While[n <= 60, $\kappa$ = Append[$\kappa$, n\^{}0 + Sum[$\kappa$[[m]], \{m, Delete[Divisors[n], -1]\}]]; n++]; $\kappa$
 }
}
\label{avalstable}
\end{table}
\subsection{Outline of paper} 
\begin{figure}[b!]
\begin{center}
\includegraphics[width=1\textwidth]{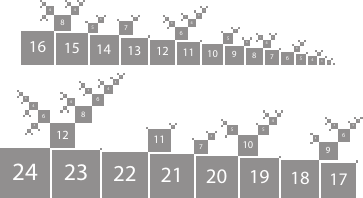}
\caption{\small
{\bf Divisor trees for 1 to 24.}
The number of recursive divisors $\kappa_0(n)$ 			counts the number of squares in each tree. 
The sum of recursive divisors $\kappa_1(n)$ 				adds up the side length of the squares in each tree. 
The sum of the square of recursive divisors $\kappa_2(n)$ 	adds up the area in each tree, and so on. 
\vspace{-14pt}
}
\label{smalltrees}
\end{center}
\end{figure}

\noindent
The goal of this paper is to quantify the notion of recursive divisibility and to understand the properties of numbers which possess it to a large degree.
It is organized as follows.
\\ \indent
In part 2, we introduce divisor trees (Figs. \ref{smalltrees} and  \ref{bigtrees}), which give a geometrical interpretation of $\kappa_x(n)$.
Using this, we find a relation between $\kappa_x(n)$ and $\kappa_0(n)$.
\\ \indent
In part 3, we show that for $x \geq 2$, $\kappa_x(n) < n^x/(2 - \zeta(x))$.
We plot $\kappa_x(n)$ for $x=1$ to $x=6$, confirming our prediction.
\\ \indent
In part 4, we investigate the number of recursive divisors $\kappa_0(n)$.
We show that $\kappa_0(n)$ is twice the number of ordered factorizations for $n \geq 2$, 
a problem much studied in its own right \cite{Kalmar,Hille,Canfield,Chor,Klazara,Deleglise}.
We give recursion relations for when $n$ is the product of distinct primes, and for when $n$ is the product of primes to a power.
The latter can be solved for up to three primes.
\\ \indent
In part 5, we investigate the sum of recursive divisors $\kappa_1(n)$.
We give recursion relations for when $n$ is the product of primes to a power.
These can be solved using the relation between $\kappa_x(n)$ and $\kappa_0(n)$ from part 2.
\\ \indent
In part 6, we study numbers which are recursively divisible to a high degree.
We call numbers with a record number of recursive divisors recursively highly composite.
These have been studied in the context of the number of ordered factorizations \cite{Deleglise}.
We call numbers with a record sum of recursive divisors, normalized by $n$, recursively super-abundant.
We list both kinds up to a million in Appendix A.
\\ \indent
In part 7, we survey applications of highly recursive numbers in design and technology and display standards.
We conclude with a list of open problems.
\section{Divisor trees and the relation between $\kappa_x$ and $\kappa_0$}  
\noindent
In this section, we prove the following relation between $\kappa_x(n)$ and $\kappa_0(n)$:
\begin{theorem}		
	\begin{equation*}
		\frac{\kappa_x(n)}{n^x} = \frac{1}{2} + \frac{1}{2} \sum_{d|n} \frac{\kappa_0(d)}{d^x}.
	\label{x0relation}
	\end{equation*}
\end{theorem} 
\noindent
To do so, we introduce the concept of divisor trees.
As well as motivating two lemmas necessary for our proof, divisor trees provide some intuition for how the recursive divisor function behaves.
\subsection{Divisor trees} 
\begin{figure}[b!]
\begin{center}
\includegraphics[width=0.8\textwidth]{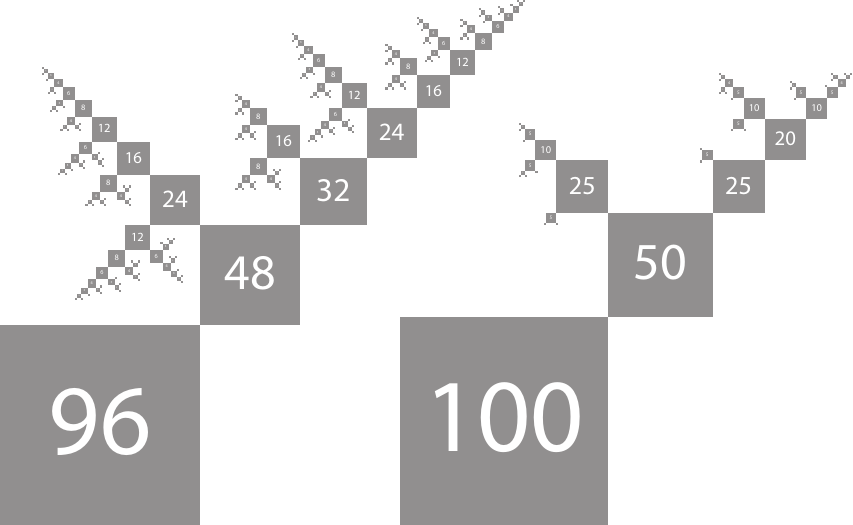}
\caption{\small
{\bf Divisor trees for 96 and 100.}
There are 		$\kappa_0(96)=224$ squares in the left tree and
		 	$\kappa_0(100)=52$ squares in the right tree.
The sum of the side length of the squares, or one-fourth of the tree perimeter, is 
			$\kappa_1(96)=768$ and
			$\kappa_1(100)=340$.
The sum of the area of the squares is 
			$\kappa_2(96)=$ 16,608 and
			$\kappa_2(100)=$ 14,740.	
\vspace{-14pt}
}
\label{bigtrees}
\end{center}
\end{figure}

A geometric interpretation of the recursive divisor function can be had by drawing the divisor tree for a given value of $n$.
Divisor trees for 1 to 24 are shown in Fig. \ref{smalltrees}. 
The number of recursive divisors $\kappa_0(n)$ counts the number of squares in each tree,
whereas the number of divisors $d \equiv \sigma_0(n)$ in (\ref{AN}) counts the number of squares in the main diagonal.
The sum of recursive divisors $\kappa_1(n)$ adds up the side length of the squares in each tree,
whereas the sum of divisors $\sigma \equiv \sigma_1(n)$ in (\ref{AN})  adds up the side length of the squares in the main diagonal.
This can be extended to $\kappa_2(n)$, which adds up area, and so on.
\\ \indent
A divisor tree is constructed as follows.
First, draw a square of side length $n$.
Let $d_1, d_2, \ldots$ be the proper divisors of $n$ in descending order. 
Then draw squares of side length $d_1$, $d_2, \ldots$ with each consecutive square situated to the upper right of its predecessor, kitty-corner, as shown in Figs. \ref{smalltrees} and \ref{bigtrees}.
This forms the main arm of a divisor tree.
Now, for each of the squares of side length $d_1, d_2, \ldots$, repeat the process.
Let $e_1, e_2, \ldots$ be the proper divisors of $d_1$ in descending order. 
Then draw squares of side length $e_1$, $e_2, \ldots$, but with the sub-arm rotated $90^{\circ}$ counter-clockwise.
Do the same for each of the remaining squares in the main arm. 
This forms the branches off of the main arm.
Now, continue repeating this process, drawing arms off of arms off of arms, and so on, until the arms are single squares of size 1.
\\ \indent
Note that, for large enough $n$, a divisor tree can overlap itself.
The precise conditions as to when is one of the open questions listed at the end.
\subsection{Properties of divisor trees} \noindent
In order to prove Theorem \ref{x0relation}, let us consider a more fine-grained description of divisor trees,
namely, one that counts the number of divisors of a given size.
\begin{definition}
	The number of recursive divisors of $n$ of size $j < n$ is
	\begin{equation*}
		\kappa_0(n,j) = 
		\begin{cases}
			\displaystyle
			\sum_{d\divprop n} \kappa_0(d,j), 	& j\divprop n\\
			0,							& {\rm otherwise},
		\end{cases}
	\end{equation*}
	where $\kappa_0(n,n)=1$ and $d\divprop n$ means $d|n$ and $d < n$.
	\label{Defank}
\end{definition}
\begin{lemma}		\label{ED}
	The number of recursive divisors of size $j$ satisifes $\kappa_0(j\,n,j) = \kappa_0(n,1)$.
\end{lemma}
\noindent
\emph{Proof.} By Definition \ref{Defank}, with $j \rightarrow 1$,
\begin{equation}
	\kappa_0(n,1) = \sum_{d\divprop n} \kappa_0(d,1).
	\label{EF}
\end{equation}
Similarly, with $n \rightarrow j \, n$,
\begin{equation*}
	\kappa_0(j\,n,j) = \sum_{d\divprop j\,n} \kappa_0(d,j).
\end{equation*}
Since $\kappa_0(d,j)=0$ if $j$ does not divide $d$, this can be rewritten as 
\begin{equation}
	\kappa_0(j\,n,j) = \sum_{d\divprop n} \kappa_0(j \, d,j).
	\label{EH}
\end{equation}
We will use this result in a moment.
\\ \indent
Let $n = p_1^{\alpha_1} \, p_2^{\alpha_2} \ldots p_\omega^{\alpha_j}$ and $\Omega(n) = \alpha_1 + \alpha_2 + \ldots + \alpha_\omega$.
We prove the lemma by induction on $\Omega(n)$. 
The base case $\Omega(n)=0$, or $n=1$, holds by Definition \ref{Defank}: $\kappa_0(j \, 1,j) = \kappa_0(1,1)$. 
We now show that if $\kappa_0(j\,n,j) = \kappa_0(n,1)$ for all $n$ such that $\Omega(n) < i$, then $\kappa_0(j\,n,j) = \kappa_0(n,1)$ for all $n$ such that $\Omega(n) < i+1$.
To see why, observe that in (\ref{EH}) all of the proper divisors $d$ of $n$ must have $\Omega(d)$ at most $\Omega(n)-1$. 
So by assumption all of the $\kappa_0(j\,d,j)$ in (\ref{EH}) reduce to $\kappa_0(d,1)$,
and the right side of (\ref{EH}) takes the form of the right side of (\ref{EF}) and thus equals $\kappa_0(n,1)$. 
\qed
\begin{lemma}		\label{FC}
For $n \geq 2$, the number of recursive divisors of size 1 is equal to half the total number of recursive divisors, that is, $\kappa_0(n,1) = \kappa_0(n)/2$.
\end{lemma}
\indent
\emph{Proof.} 
Clearly
\begin{equation*}
\kappa_0(n) = \sum_{d|n} \kappa_0(n,d).
\end{equation*}
By Lemma \ref{ED}, $\kappa_0(n,k) = \kappa_0(n/k,1)$ for $k|n$, so the above becomes
\begin{align*}
	\kappa_0(n) 	&= \sum_{d|n} \kappa_0(n/d,1) \\
				&= \sum_{d|n} \kappa_0(d,1) \\
				&= \kappa_0(n,1) + \sum_{d\divprop n} \kappa_0(d,1).
\end{align*}
Inserting Definition \ref{Defank} with $j=1$ into the above gives the desired result. \qed
\subsection{Relation between \boldmath{$\kappa_x$} and $\kappa_0$.} 
\noindent 
\emph{Proof of Theorem \ref{x0relation}}.
We can write $\kappa_x(n)$ as
\begin{equation*}
\kappa_x(n) = \sum_{d|n} d^x \, \kappa_0(n,d).
\label{FB}
\end{equation*}
By Lemma \ref{ED},
\begin{equation*}
\kappa_x(n) = \sum_{d|n} d^x \, \kappa_0(n/d,1).
\label{FD}
\end{equation*}
Recall that Lemma \ref{FC} only applies for $n \geq 2$, so we pull out the $d = n$ term:
\begin{equation*}
\kappa_x(n) = n^x + \sum_{d\divprop n} d^x \, \kappa_0(n/d,1).
\end{equation*}
By Lemma \ref{FC},
\begin{align*}
\kappa_x(n) 	&= n^x + \frac{1}{2} \sum_{d\divprop n} d^x \, \kappa_0(n/d) \\
			&= \frac{n^x}{2} + \frac{1}{2} \sum_{d|n} d^x \kappa_0(n/d) \\
			&= \frac{n^x}{2} + \frac{1}{2} \sum_{d|n} \left(\frac{n}{d}\right)^x \kappa_0(d).
\end{align*}
Dividing by $n^x$, the theorem follows. \qed
\begin{figure}[b!]
\includegraphics[width=1\textwidth]{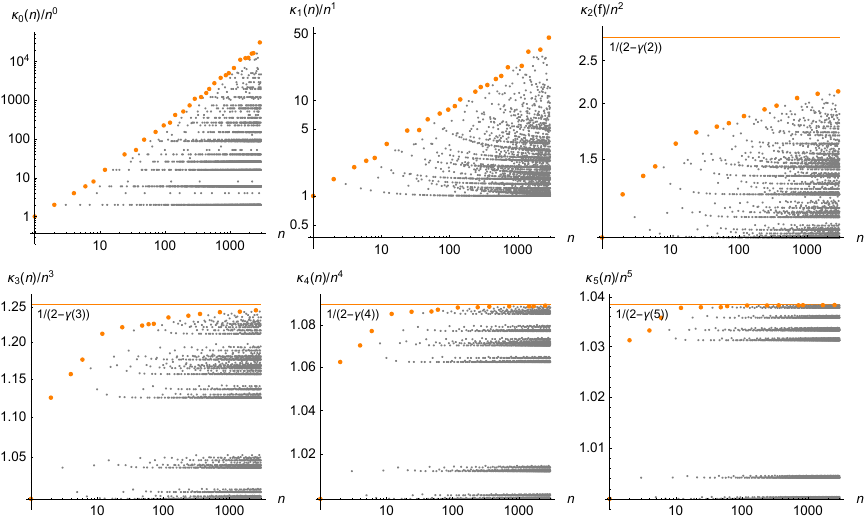} 
\caption{\small
{\bf Plots of} 
\boldmath
$\kappa_x(n)/n^x.$
\unboldmath
We plot $\kappa_x(n)/n^x$ up to $n =$ 3,000, for $x=0$ to $x=5$.
The large orange points are the sequence records, which satisfy $\kappa_x(n)/n^x > \kappa_x(m)/m^x$ for all $m < n$.
For $x=0$ and $x=1$, the upper bound for the sequence diverges.
But for $x=2$ and above, it converges to $1/(2-\gamma(x))$, where $\gamma$ is the Riemann zeta function.
For $x=0$, the orange points are the recursively highly composite numbers.
For $x=1$, they are the recursively super-abundnant numbers (see Appendix A).
}
\vspace{-14pt}
\label{avalsplot}
\end{figure}
\section{Properties of $\kappa_2$, $\kappa_3$, and so on.}
\noindent
\begin{theorem}
	For $x > 1$, 
		\begin{equation*}
		\frac{\kappa_x(n)}{n^x} < \frac{1}{2 - \zeta(x)}.
		\end{equation*}
	\end{theorem}
\noindent
This theorem was inspired by an anonymous referee from a previous version of this paper.
The referee answered one of our then open questions at the end of our paper.
We generalized the answer, which led to this theorem.
We confirm it in Fig. \ref{avalsplot}, in which we plot $\kappa_x(n)/n^x$ for $x=1$ to $x=6$. 
\\ \indent
\emph{Proof.}
We prove the theorem by induction.
We know that it is true for $n=1$.
Assume it is true for all numbers less than $n$.
We show it is then true for $n$.
By Definition \ref{MainDef}, 
\begin{eqnarray*}
\frac{\kappa_x(n)}{n^x} 	&=&		1 + \frac{1}{n^x} \sum_{d\divprop n} \kappa_x(d)  \\
					&=&		1 + \sum_{d\divprop n} \frac{\kappa_x(d)}{d^x} \frac{d^x}{n^x} \\
					&<&		1 + \frac{1}{2 - \zeta(x)} \sum_{d\divprop n} \frac{d^x}{n^x} \\
					&<&		1 + \frac{1}{2 - \zeta(x)} \Bigg(-1 + \sum_{d\vert n} \frac{d^x}{n^x} \Bigg),
\end{eqnarray*}
where $\gamma(x)$ is the Riemann zeta function.
Since, for $x > 1$, 
\begin{equation}
\sum_{d\vert n} \frac{d^x}{n^x} =  \sum_{d\vert n} \frac{1}{d^x} <  \zeta(x),
\end{equation}
we have
\begin{eqnarray*}
\frac{\kappa_x(n)}{n^x} 	&<&		1 + \frac{1}{2 - \zeta(x)} \left(\zeta(x) - 1\right) \\
					&<&		\frac{1}{2 - \zeta(x)}. \quad \qed
\end{eqnarray*}
\section{Number of recursive divisors}
\subsection{Relation to ordered factorizations} 
\noindent
Some of the properties of $\kappa_0(n)$ can be deduced from properties of a related function, the number $K(n)$ of ordered factorizations into integers greater than one.
It satisfies $K(1)=1$ and
	\begin{equation*}
		K(n) = \sum_{d\divprop n} K(d).
		\label{AL}
	\end{equation*} 
For example, 12 is the product of integers greater than one in eight ways: $12=6 \cdot 2=2 \cdot 6=4 \cdot 3=3 \cdot 4=3 \cdot 2 \cdot 2=2 \cdot 3 \cdot 2=2 \cdot 2 \cdot 3$.
So $K(12)=8$.
Kalmar \cite{Kalmar} was the first to consider $K(n)$---hence the name $K$---and it was later studied more systematically by Hille \cite{Hille}. 
Over the last 20 years several authors have extended Hille's results \cite{Chor,Klazara,Deleglise}, some of which we will mention later.
\\ \indent
Notice that the definition of $K(n)$ is identical to Definition \ref{Defank} for $j=1$, that is, to $\kappa_0(n,1)$.
Since $K(1) = \kappa_0(1,1) = 1$, then by Lemma \ref{FC} we arrive at
	\begin{proposition}
		For $n\geq2$, $\kappa_0(n) = 2 K(n)$, where $K(n)$ is the number of ordered factorizations into integers greater than one.
		\label{avsg}
	\end{proposition}
\subsection{Distinct primes} \noindent
Let $n = p_1 p_2 \ldots p_\omega$ be the product of $\omega$ distinct primes. 
\begin{proposition}
	The exponential generating function of $\kappa_0(p_1 \ldots p_\omega)$ is
	\begin{equation*}
		{\rm EG}(\kappa_0(p_1 \ldots p_\omega),x)= 	 \frac{e^x}{2-e^x}.
	\end{equation*}
\end{proposition}
\emph{Proof.}
Since $\kappa_0$ depends only on the prime signature of $n$, which in this case is all 1s, we can immediately write down
	\begin{equation*}
		2 \kappa_0(p_1 \ldots p_\omega)	=  1 + \sum_{i=0}^{\omega} \binom{\omega}{i} \kappa_0(p_1 \ldots p_i).
	\end{equation*}
Then
	\begin{align*}
		2 \sum_{\omega=0}^\infty \frac{x^w}{w!} \kappa_0(p_1,\ldots,p_\omega) 
		&=  \sum_{\omega=0}^\infty \frac{x^\omega}{\omega!} \left(	1 + \sum_{i=0}^{\omega} \binom{\omega}{i} \kappa_0(p_1 \ldots p_i) \right) \\
		&=  e^x + \sum_{\omega=0}^\infty  \sum_{i=0}^{\omega} \frac{x^i}{i!} \frac{x^{\omega-i}}{(\omega-i)!} \kappa_0(p_1 \ldots p_i).
	\end{align*}
Summing this along the diagonals $\omega = i$, $\omega = i + 1$, and so on,
	\begin{align*}
		2 \sum_{\omega=0}^\infty \frac{x^w}{w!} \kappa_0(p_1,\ldots,p_\omega) 
		&\!=\! e^x 
		\!+\! \frac{x^0}{0!} \sum_{\omega=0}^\infty \frac{x^w}{w!}  \kappa_0(p_1,\ldots,p_\omega)
		\!+\! \frac{x^1}{1!} \sum_{\omega=0}^\infty \frac{x^w}{w!}  \kappa_0(p_1,\ldots,p_\omega) 
		\!+\! \ldots \\
		&=  e^x + e^x \sum_{\omega=0}^\infty \frac{x^w}{w!}  \kappa_0(p_1, \ldots, p_\omega).
	\end{align*}
From this it follows that 
	\begin{equation*}
		{\rm EG}(\kappa_0(p_1 \ldots p_\omega),x)= 	 \frac{e^x}{2-e^x}. \qed
	\end{equation*}
\subsection{Primes to a power} \noindent
When $n$ is the product of primes to powers, $\kappa_0(n)$ satisfies recursion relations relating it to values of $\kappa_0(n)$ for primes to lower powers.
The first few of these can be solved explicitly.
\begin{theorem}
Let $p, q$ and $r$ be prime. Then 
	\begin{align*}
	\kappa_0(p^c) 		&= \textstyle 2 \kappa_0\big(\frac{p^c}{p}\big) \\ &= 2^c																		\\
	\kappa_0(p^c q^d) 	&= 2 \left( \textstyle \kappa_0\big(\frac{p^c q^d}{p}\big) + \kappa_0\big(\frac{p^c q^d}{q}\big) - \kappa_0\big(\frac{p^c q^d}{pq}\big) \right) 		\\
					&= 2^c \sum_{i=0}^d \textstyle \binom{d}{i} \binom{c+i}{i},	 \\
	\kappa_0(p^c q^d r^e) &= 2 \Big(\textstyle \kappa_0\big(\frac{p^c q^d r^e}{p}\big) + \kappa_0\big(\frac{p^c q^d r^e}{q}\big) + \kappa_0\big(\frac{p^c q^d r^e}{r}\big) \\
					 	& \textstyle - \kappa_0\big(\frac{p^c q^d r^e}{pq}\big) - \kappa_0\big(\frac{p^c q^d r^e}{pr}\big) - \kappa_0\big(\frac{p^c q^d r^e}{qr}\big) + \kappa_0\big(\frac{p^c q^d r^e}{pqr}\big) \Big) \\
					&= 		\sum_{j=0}^d (-1)^j \textstyle \binom{d}{j} \binom{c+d-j}{d} \kappa_0(p^{c+d-j} r^e).
	\end{align*}
Analogous recursion relations apply for the product of more primes to powers.
	\label{arec}
\end{theorem}
\indent
\emph{Proof.}
For the recursion relations, the approach is similar to, but somewhat simpler than, that used to prove the recursion relations in Theorem \ref{brec}.
But we can deduce these relations from previous work using Proposition \ref{avsg}.
Hille \cite{Hille} and Chor \emph{et al.} \cite{Chor} proved that identical recursion relations govern $K(n)$, the number of ordered factorizations. 
Since $\kappa_0(n) = 2 K(n)$ for $n \geq 2$, the same recursion relations apply to $\kappa_0(n)$.
For the explicit forms of $\kappa_0$, Chor \emph{et al.} \cite{Chor} give analogous results for $K(n)$, which when multiplied by 2 apply to $\kappa_0(n)$.
\qed
\begin{corollary}
Let $\alpha^*$ be the maximum exponent in the prime factorization of $n$. 
Then $2^{\alpha^*}$ divides $\kappa_0(n)$.
\end{corollary}
\emph{Proof.}
All of the recursion relations in Theorem \ref{arec} have a factor of 2 on the right side.
The corollary is implied by iterating the recursion relation $\alpha^*$ times.
Each time, the exponents on the right are reduced by at most 1.
Iterating until the smallest exponent is reduced to 0, the exponent disappears since, for example, $\kappa_0(p^c q^0) = \kappa_0(p^c)$.
Continuing this process ultimately gives a total of $\alpha^*$ factors of 2. 
The $\kappa_0(n)$ are expressed as a product of an integer and $2^{\alpha^*}$ in Appendix A.
\qed
\section{Sum of recursive divisors}
\noindent 
We now turn to the sum of recursive divisors $\kappa_1(n)$.
This quantity is more intricate than $\kappa_0(n)$, because it depends on the primes as well as their exponents  in the prime factorization of $n$.
\subsection{Primes to a power} \noindent
When $n$ is equal to the product of primes to powers, $\kappa_1(n)$ satisfies recursion relations similar to those for $\kappa_0(n)$, but more complex.
\begin{theorem}		\label{brec}
Let $p, q$ and $r$ be prime. Then
	\begin{align*}
	\kappa_1(p^c) 		& = \textstyle 2 \kappa_1\big(\frac{p^c}{p}\big) + \frac{p-1}{p} p^c \\ 
					& = \begin{cases}
						2^c \, \frac{c+2}{2} 			& p=2 \\
						p^c \, \frac{p-1 - (2/p)^c}{p-2} 	& p \; {\rm odd}
						\end{cases}																								\\
\kappa_1(p^c q^d) 		&=  \textstyle 2 \left(\kappa_1\big(\frac{p^c q^d}{p}\big) + \kappa_1\big(\frac{p^c q^d}{q}\big) - \kappa_1\big(\frac{p^c q^d}{pq}\big) \right) +  \frac{p-1}{p} \frac{q-1}{q} p^c q^d \\
					&= p^c q^d \bigg(\frac{1}{2} + \frac{1}{2} \sum_{i=0}^c \frac{2^i}{p^i} \sum_{j=0}^d \frac{1}{q^j} \sum_{k=0}^j \textstyle \binom{i+k}{k} \binom{j}{k} \bigg)	 \\
\kappa_1(p^c q^d r^e) 	&= 2 \Big(\textstyle \kappa_1\big(\frac{p^c q^d r^e}{p}\big) + \kappa_1\big(\frac{p^c q^d r^e}{q}\big) + \kappa_1\big(\frac{p^c q^d r^e}{r}\big) \\
					& \textstyle - \kappa_1\big(\frac{p^c q^d r^e}{pq}\big) - \kappa_1\big(\frac{p^c q^d r^e}{pr}\big) - \kappa_1\big(\frac{p^c q^d r^e}{qr}\big) 
					+ \kappa_1\big(\frac{p^c q^d r^e}{pqr}\big) \Big) \\
					& \textstyle + \frac{p-1}{p} \frac{q-1}{q} \frac{r-1}{r} p^c q^d r^e.
	\end{align*}
\end{theorem}
\emph{Proof.}
We start with the recursion relations.
For $n=p^c$, from Definition \ref{MainDef}, 
\begin{equation}
\kappa_1(p^c) = p^c + \sum_{i=0}^{c-1} \kappa_1(p^i).
\label{SB}
\end{equation} 
Adding $\kappa_1(p^c)$ to both sides and with $c \rightarrow c-1$, 
\begin{equation*}
\sum_{i=0}^{c-1} \kappa_1(p^i) = 2 \kappa_1(p^{c-1}) - p^{c-1},
\label{SF}
\end{equation*} 
which, when inserted into (\ref{SB}), gives the desired recurrence relation.
\\ \indent
For $n=p^c q^d$, from Definition \ref{MainDef}, 
\begin{equation}
\kappa_1(p^{c} q^{d}) = p^c q^d + \sum_{i=0}^{c-1} \sum_{j=0}^{d} \kappa_1(p^{i} q^{j}) + \sum_{j=0}^{d-1} \kappa_1(p^{c} q^j).
\label{TB}
\end{equation} 
Adding $\kappa_1(p^{c} q^{d})$ to both sides,
\begin{equation}
2 \kappa_1(p^c q^d) = p^c q^d + \sum_{i=0}^{c-1} \sum_{j=0}^d \kappa_1(p^{i} q^{j}) + \sum_{j=0}^d \kappa_1(p^c q^j),
\label{TD}
\end{equation} 
which we can equally write
\begin{equation}
2 \kappa_1(p^c q^d) = p^c q^d + \sum_{i=0}^{c} \sum_{j=0}^d \kappa_1(p^i q^j).
\label{TF}
\end{equation}
With $d \rightarrow d-1$ in (\ref{TD}), we find
\begin{equation}
 \sum_{j=0}^{d-1} \kappa_1(p^c q^j) = 2 \kappa_1(p^{c} q^{d-1}) - p^c q^{d-1} - \sum_{j=0}^{c-1} \sum_{i=0}^{d-1} \kappa_1(p^i q^j).
 \label{TH}
\end{equation} 
With $c \rightarrow c-1$ and $d \rightarrow d-1$ in (\ref{TF}), and inserting the result into (\ref{TH}), yields 
\begin{equation}
 \sum_{j=0}^{d-1} \kappa_1(p^c q^j) = 2 \kappa_1(p^{c} q^{d-1}) - 2 \kappa_1(p^{c-1} q^{d-1}) +(1-p) p^{c-1} q^{d-1}.
 \label{TL}
\end{equation}
With $c \rightarrow c-1$ in (\ref{TF}), we find
\begin{equation}
 \sum_{i=0}^{c-1} \sum_{j=0}^{d} \kappa_1(p^i q^j) = 2 \kappa_1(p^{c-1} q^{d}) - p^{c-1} q^d.
\label{TN}
\end{equation}
Inserting (\ref{TL}) and (\ref{TN}) into (\ref{TB}) gives the desired recursion relation.
\\ \indent
For $n = p^c q^d r^e$, the proof is similar to the one above and is omitted here.
\\ \indent
For the explicit forms of $\kappa_1(n)$, we appeal to Theorem \ref{x0relation}, which tells us
\begin{equation}
	\frac{2 \kappa_1(n)}{n} = 1 + \sum_{d|n} \frac{\kappa_0(d)}{d}.
	\label{VB}
\end{equation}
For $n=p^c$, from (\ref{VB}) we have
	\begin{equation*}
		\frac{2 \kappa_1(p^c)}{p^c} = 1 + \sum_{i=0}^c \frac{\kappa_0(p^i)}{p^i}.
	\end{equation*}
Inserting Theorem \ref{arec} into the above gives the desired result.
For $n=p^c q^d$, from (\ref{VB}) we have
	\begin{equation*}
		\frac{2 \kappa_1(p^c q^d)}{p^c q^d} = 1 + \sum_{i=0}^c \sum_{j=0}^d \frac{\kappa_0(p^i q^j)}{p^i q^j}.
	\end{equation*}
Inserting Theorem \ref{arec} into the above gives the desired result.
For $p=2$, the result simplifies to contain just two sums. \qed
\section{Numbers that are recursively divisible to a high degree} \noindent
\subsection{Recursively highly composite numbers} \noindent
A number $n$ is highly composite \cite{RamanujanA} if it has more divisors than any of its predecessors, that is, $\sigma_0(n) > \sigma_0(m)$ for all $m < n$.
These are shown in the right side of Table 2 in Appendix A.
\\ \indent
By analogy with highly composite numbers, a number $n$ is recursively highly composite if it has more recursive divisors than any of its predecessors.
\begin{definition}
$n$ is recursively highly composite if $\kappa_0(n) > \kappa_0(m)$ for all $m < n$.
\end{definition}
\noindent
These numbers are shown in the left side of Table 2 in Appendix A.
From the third term, they correspond to the indices of sequence records of $K(n)$, the K-champion numbers \cite{Deleglise}.
Because $\kappa_0(n)$ depends only on the exponents in the prime factorization of $n$, the exponents in recursively highly composite numbers must be non-increasing.
\subsection{Recursively super-abundant numbers} \noindent
A number $n$ is super-abundant \cite{ErdosA} if the sum of its divisors, normalized by $n$, is greater than that of any of its predecessors, that is, $\sigma(n)/n > \sigma(m)/m$ for all $m < n$.
These are the starred numbers in the right side of Table 2 in Appendix A.
For small $n$, super-abundant numbers are also highly composite, but later this ceases to be the case.
The first super-abundant number that is not highly composite is 1,163,962,800 (A166735 \cite{Sloane}), and in fact only 449 numbers have both properties (A166981 \cite{Sloane}).
\\ \indent
By analogy with super-abundant numbers, a number $n$ is recursively super-abundant if the sum of its recursive divisors, normalized by $n$, is greater than that of any of its predecessors.
\begin{definition}
$n$ is recursively super-abundant if $\kappa_1(n)/n > \kappa_1(m)/m$ for all $m < n$.
\end{definition}
\noindent
These numbers are starred in the left side of Table 2 in Appendix A.
Early on, recursively super-abundant numbers are recursively highly composite. The first exception is 181,440.
\setcounter{table}{2}
\setlength{\tabcolsep}{5pt}
\begin{table*}[b!]
\begin{small}
\begin{tabular}{rrlrl}
$n$ 		& \multicolumn{2}{c}{\emph{Design and technology}} 		& \multicolumn{2}{c}{\emph{Display standards}} 		\vspace{1pt}	\\
*24   		& $24 \times 16$ 			& Biotech 384-well assay 		& 						&				\vspace{0pt} 	\\
*48  		& $128 \times 48$			& TRS 80					& 						&				\vspace{0pt} 	\\
72 		& 72 points/in				& Adobe typography point		& 						&				\vspace{0pt} 	\\
96		& $96 \times 65$			& Nokia 1100 phone			& 				 		&				\vspace{0pt}	\\
*120	 	& $120 \times 160$			& Nokia 100 phone				& $160 \times 120$			& QQVGA 	\vspace{0pt} 	\\
144 		& $144 \times 168$			& Pebble Time watch		& 						&				\vspace{0pt} 	\\
*240 		& $240 \times 64$			& Atari Portfolio				& $320 \times 240$			& Quarter VGA  	\vspace{0pt} 	\\
288 		& $352 \times 288$			& Video CD 				& $352 \times 288$			& CIF			\vspace{0pt} 	\\
*360 		& $360 \times 360$			& LG Watch Style 			& $640 \times 360$			& nHD 			\vspace{0pt} 	\\
480 		& $320 \times 480$			& iPhone 1--3				& $640 \times 480$			& VGA			\vspace{0pt} 	\\
576 		& 576 lines				& PAL analog television		& $1024 \times 576$			& WSVGA 		\vspace{0pt} 	\\
*720 		& $720 \times 364$			& Macintosh XL, Hercules 	& $1280 \times 720$ 		& HD  			\vspace{0pt} 	\\
864		& 						&						& $1152 \times 864$			& XGA+			\vspace{0pt} 	\\		
960 		&						& Facebook website 			& 						&				\vspace{0pt} 	\\
*1152	&						& 						& $1152 \times 2048$		& QWXGA 		\vspace{0pt} 	\\
*1440 	&						& 3.5" disk block size		& $2560 \times 1440$		& Quad HD 		\vspace{0pt} 	\\
1920		&						&						& $1920 \times 1080$		& Full HD			\vspace{0pt} 	\\
*2160 	& $2160 \times 1440$ 		& Microsoft Surface Pro 3		& $4096 \times 2160$		& 4K Ultra HD		\vspace{0pt} 	\\
2304 	& $2304 \times 1440$ 		& MacBook Retina 			& 						&				\vspace{0pt} 	\\
*2880 	& $\quad 2880 \times 1800$	& 15" MacBook Pro Retina 	& $\quad$ $5120 \times 2880$	& 5K  			\vspace{0pt} 	\\
3456 	&						& Canon EOS 1100D		& 						&				\vspace{0pt} 	\\
*4320 	&						& 						& $7680 \times 4320$		& 8K Ultra HD  		\vspace{0pt} 	\\
*8640 	&						& 						& $15360 \times 8640$		&16K Ultra HD 		\vspace{0pt} 	\\
\end{tabular}
\vspace{0.1in}
\caption{\small 
{\bf Applications.}
Numbers that are recursively divisible to a large degree predict the numbers that frequently show up in design and technology and display standards.
All of the numbers $n$ are recursively highly composite; those that are starred are also recursively super-abundant.
}
\end{small}
\label{applications}
\label{Applications}
\end{table*}
\subsection{Applications.} \noindent
In graphic and digital design, the layout of graphics and text is often constrained to lie on an underlying rectangular grid \cite{Brockmann}.
The grid elements are the primitive building blocks from which bigger columns or rows can be formed.
For example, grids of 24 and 96 columns are often used for books and websites, respectively \cite{Brockmann}.
Using a grid reduces the space of possible designs, making it easier to navigate.
And the design elements become more interoperable, like how Lego bricks snap into place with one another, making it faster to build new designs. 
\\ \indent
What are the best grid sizes?
The challenge is committing to a grid size now that provides the greatest optionality for an unknown future.
Imagine, for example, that we have to cut a pie into slices, to be divided up later for an unknown number of colleagues.
How many slices should we choose?
The answer in this case is straightforward: the best grids are the ones with the most divisors, such as the highly composite or super-abundant numbers \cite{RamanujanA,ErdosA}.
\\ \indent
But the story gets more complicated when we need to consider steps into the future.
For instance, imagine now that each colleague takes his share of pie home to further divide it amongst his family---but they cannot make any additional cuts. 
In this case, not only does the whole need to be highly divisible, but the parts need to be highly divisible, too.
This process can be extended in a recursive way.
\\ \indent
Recursive modularity, in which there are multiple levels of organization, has long been a feature of graphic and digital design.
For example, newspapers are divided into columns for different stories, and columns into sub-columns of text.
But with the rise of digital technologies, recursive modularity is becoming the rule.
Different pages of a website are divided into different numbers of columns, each of which can be broken down into smaller design elements.
Often one column from the website fills the full screen of a phone. 
\\ \indent
Specific applications of recursively highly composite numbers are shown in Table 3.
In design and technology, these numbers are used for the screen resolutions of watches, phones, cameras and computers.
They appear in typesetting, websites, and experimental equipment, such as test tube microplates.
In display standards, many resolutions use these numbers in the height or width, measured in pixels. 
Because these standards tend to preserve certain aspect ratios, such as 16:9, usually just one of the two dimensions is highly recursively divisible.
\subsection{Open questions} \noindent
There are many open questions about the recursive divisor function and numbers that are recursively divisible to a high degree.
Here are six.
\\ 1. Let $n$ be a number such that $\sigma_1(n) > 3 n$. 
The first such numbers are 180, 240, 360, \ldots (A068403 \cite{Sloane}).
Then the divisor tree for $4 n$ overlaps itself (see lims.ac.uk/recursively-divisible-numbers).
But there are other numbers that cause overlaps. 
What are they?
\\ 2. For what values of $n$ do divisor trees have an (approximate) fractal dimension? 
\\ 3. What is the value of $\kappa_1(n)$ when $n$ is the product of the first $k$ distinct primes?
\\ 4. How frequently do recursively highly composite numbers appear? How about recursively super-abundant numbers? 
\\ 5. Recursively perfect numbers satisfy $\kappa_0(n) = n$ \cite{Fink22}. How frequently do they appear?
\\ 6. Recursively abundant numbers satisfy $\kappa_0(n) > n$ \cite{Fink22}. Are any odd and, if so, what is the smallest?
\\ \indent
I acknowledge Andriy Fedosyeyev for creating the divisor tree generator, \\ lims.ac.uk/recursively-divisible-numbers.
\vspace{1.6in}
\section*{Appendix A} 
\vspace{-3mm}
\begin{small}
\begin{center}
$
\begin{array}{rclrcrclr}
	&	n	&		&.    \mbox{\hspace{0.4in}} \kappa_0(n)	& 	\mbox{\hspace{0.5in}} 	&		&   n  	& 		&    \mbox{\hspace{0.2in}} \sigma_0(n) \\
 *1 & = & 1 &1 & \text{} 									& *1 & = & 1 & 1 \\
 *2 & = & 2 & 1 \cdot 2 & \text{} 								& *2 & = & 2 & 2 \\
 *4 & = & 2^2 & 1\cdot 2^2 & \text{} 							& *4 & = & 2^2 & 3 \\
 *6 & = & 2\cdot 3 & 3\cdot 2 & \text{} 						& *6 & = & 2\cdot 3 & 4 \\
 8 & = & 2^3 & 1\cdot 2^3 & \text{} 							& \text{} & \text{} & \text{} & \text{} \\
 *12 & = & 2^2\cdot 3 & 4\cdot 2^2 & \text{} 					& *12 & = & 2^2\cdot 3 & 6 \\
 *24 & = & 2^3\cdot 3 & 5\cdot 2^3 & \text{} 					& *24 & = & 2^3\cdot 3 & 8 \\
 *36 & = & 2^2\cdot 3^2 & 13\cdot 2^2 & \text{} 					& *36 & = & 2^2\cdot 3^2 & 9 \\
 *48 & = & 2^4\cdot 3 & 6\cdot 2^4 & \text{} 					& *48 & = & 2^4\cdot 3 & 10 \\
 \text{} & \text{} & \text{} & \text{} & \text{} 						& *60 & = & 2^2\cdot 3\cdot 5 & 12 \\
 72 & = & 2^3\cdot 3^2 & 19\cdot 2^3 & \text{} 					& \text{} & \text{} & \text{} & \text{} \\
 96 & = & 2^5\cdot 3 & 7\cdot 2^5 & \text{} 					& \text{} & \text{} & \text{} & \text{} \\
 *120 & = & 2^3\cdot 3\cdot 5 & 33\cdot 2^3 & \text{} 			& *120 & = & 2^3\cdot 3\cdot 5 & 16 \\
 144 & = & 2^4\cdot 3^2 & 26\cdot 2^4 & \text{} 					& \text{} & \text{} & \text{} & \text{} \\
 \text{} & \text{} & \text{} & \text{} & \text{} 						& *180 & = & 2^2\cdot 3^2\cdot 5 & 18 \\
 192 & = & 2^6\cdot 3 & 8\cdot 2^6 & \text{} 					& \text{} & \text{} & \text{} & \text{} \\
 *240 & = & 2^4\cdot 3\cdot 5 & 46\cdot 2^4 & \text{} 			& *240 & = & 2^4\cdot 3\cdot 5 & 20 \\
 288 & = & 2^5\cdot 3^2 & 34\cdot 2^5 & \text{} 					& \text{} & \text{} & \text{} & \text{} \\
 *360 & = & 2^3\cdot 3^2\cdot 5 & 151\cdot 2^3 & \text{} 			& *360 & = & 2^3\cdot 3^2\cdot 5 & 24 \\
 432 & = & 2^4\cdot 3^3 & 96\cdot 2^4 & \text{} 					& \text{} & \text{} & \text{} & \text{} \\
 480 & = & 2^5\cdot 3\cdot 5 & 61\cdot 2^5 & \text{} 				& \text{} & \text{} & \text{} & \text{} \\
 576 & = & 2^6\cdot 3^2 & 43\cdot 2^6 & \text{} 					& \text{} & \text{} & \text{} & \text{} \\
 *720 & = & 2^4\cdot 3^2\cdot 5 & 236\cdot 2^4 & \text{} 			& *720 & = & 2^4\cdot 3^2\cdot 5 & 30 \\
 \text{} & \text{} & \text{} & \text{} & \text{} 						& *840 & = & 2^3\cdot 3\cdot 5\cdot 7 & 32 \\
 864 & = & 2^5\cdot 3^3 & 138\cdot 2^5 & \text{} 				& \text{} & \text{} & \text{} & \text{} \\
 960 & = & 2^6\cdot 3\cdot 5 & 78\cdot 2^6 & \text{} 				& \text{} & \text{} & \text{} & \text{} \\
 *1152 & = & 2^7\cdot 3^2 & 53\cdot 2^7 & \text{} 				& \text{} & \text{} & \text{} & \text{} \\
 \text{} & \text{} & \text{} & \text{} & \text{} 						& *1260 & = & 2^2\cdot 3^2\cdot 5\cdot 7 & 36 \\
 *1440 & = & 2^5\cdot 3^2\cdot 5 & 346\cdot 2^5 & \text{} 			& \text{} & \text{} & \text{} & \text{} \\
 \text{} & \text{} & \text{} & \text{} & \text{} 						& *1680 & = & 2^4\cdot 3\cdot 5\cdot 7 & 40 \\
 1728 & = & 2^6\cdot 3^3 & 190\cdot 2^6 & \text{} 				& \text{} & \text{} & \text{} & \text{} \\
 1920 & = & 2^7\cdot 3\cdot 5 & 97\cdot 2^7 & \text{} 			& \text{} & \text{} & \text{} & \text{} \\
 *2160 & = & 2^4\cdot 3^3\cdot 5 & 996\cdot 2^4 & \text{} 			& \text{} & \text{} & \text{} & \text{} \\
 2304 & = & 2^8\cdot 3^2 & 64\cdot 2^8 & \text{} & \text{} 			& \text{} & \text{} & \text{} \\
 \text{} & \text{} & \text{} & \text{} & \text{} 						& *2520 & = & 2^3\cdot 3^2\cdot 5\cdot 7 & 48 \\
  *2880 & = & 2^6\cdot 3^2\cdot 5 & 484\cdot 2^6 & \text{} 		& \text{} & \text{} & \text{} & \text{} \\
   3456 & = & 2^7\cdot 3^3 & 253\cdot 2^7 & \text{} 				& \text{} & \text{} & \text{} & \text{} \\
    *4320 & = & 2^5\cdot 3^3\cdot 5 & 1590\cdot 2^5 & \text{} 		& \text{} & \text{} & \text{} & \text{} \\
 \text{} & \text{} & \text{} & \text{} & \text{} 						& *5040 & = & 2^4\cdot 3^2\cdot 5\cdot 7 & 60 \\
 *5760 & = & 2^7\cdot 3^2\cdot 5 & 653\cdot 2^7 & \text{} 			& \text{} & \text{} & \text{} & \text{} \\
 6912 & = & 2^8\cdot 3^3 & 328\cdot 2^8 & \text{} 				& \text{} & \text{} & \text{} & \text{} \\
 \text{} & \text{} & \text{} & \text{} & \text{} 						& 7560 & = & 2^3\cdot 3^3\cdot 5\cdot 7 & 64 \\
  *8640 & = & 2^6\cdot 3^3\cdot 5 & 2402\cdot 2^6 & \text{} 		& \text{} & \text{} & \text{} & \text{} \\
 \text{} & \text{} & \text{} & \text{} & \text{} 						& *10080 & = & 2^5\cdot 3^2\cdot 5\cdot 7 & 72 
 \end{array}
$
 \end{center}
TABLE 2. 
{\bf Numbers that are recursively divisible to a high degree.} 
The left side shows 		the recursively highly composite numbers and 
					the recursively super-abundant numbers (starred)
					up to a million.
					All of the recursively super-abundant numbers shown are also recursively highly composite, apart from one, 181,440.
The right side shows 	the highly composite numbers and 
					the super-abundant numbers (starred)
					up to a million.
					All of the super-abundant numbers shown are also highly composite.
\end{small}
\begin{small}
\begin{center}
$
\begin{array}{rclrcrclr}
	&	n	&		&    \mbox{\hspace{0.1in}}  \kappa_0(n)	& 	\mbox{\hspace{0.0in}} 	&		&   n  	& 		&    \mbox{\hspace{0.0in}} \sigma_0(n) \\
  *11520 & = & 2^8\cdot 3^2\cdot 5 & 856\cdot 2^8 & \text{} 		& \text{} & \text{} & \text{} & \text{} \\
 \text{} & \text{} & \text{} & \text{} & \text{} & *15120 & = & 2^4\cdot 3^3\cdot 5\cdot 7 & 80 \\
 *17280 & = & 2^7\cdot 3^3\cdot 5 & 3477\cdot 2^7 & \text{} & \text{} & \text{} & \text{} & \text{} \\
 \text{} & \text{} & \text{} & \text{} & \text{} & 20160 & = & 2^6\cdot 3^2\cdot 5\cdot 7 & 84 \\
 23040 & = & 2^9\cdot 3^2\cdot 5 & 1096\cdot 2^9 & \text{} & \text{} & \text{} & \text{} & \text{} \\
 \text{} & \text{} & \text{} & \text{} & \text{} & *25200 & = & 2^4\cdot 3^2\cdot 5^2\cdot 7 & 90 \\
 *25920 & = & 2^6\cdot 3^4\cdot 5 & 10368\cdot 2^6 & \text{} & \text{} & \text{} & \text{} & \text{} \\
 \text{} & \text{} & \text{} & \text{} & \text{} & *27720 & = & 2^3\cdot 3^2\cdot 5\cdot 7\cdot 11 & 96 \\
 *30240 & = & 2^5\cdot 3^3\cdot 5\cdot 7 & 20874\cdot 2^5 & \text{} & \text{} & \text{} & \text{} & \text{} \\
 *34560 & = & 2^8\cdot 3^3\cdot 5 & 4864\cdot 2^8 & \text{} & \text{} & \text{} & \text{} & \text{} \\
 \text{} & \text{} & \text{} & \text{} & \text{} & 45360 & = & 2^4\cdot 3^4\cdot 5\cdot 7 & 100 \\
 46080 & = & 2^{10}\cdot 3^2\cdot 5 & 1376\cdot 2^{10} & \text{} & \text{} & \text{} & \text{} & \text{} \\
 \text{} & \text{} & \text{} & \text{} & \text{} & 50400 & = & 2^5\cdot 3^2\cdot 5^2\cdot 7 & 108 \\
 *51840 & = & 2^7\cdot 3^4\cdot 5 & 15979\cdot 2^7 & \text{} & \text{} & \text{} & \text{} & \text{} \\
 \text{} & \text{} & \text{} & \text{} & \text{} & *55440 & = & 2^4\cdot 3^2\cdot 5\cdot 7\cdot 11 & 120 \\
 *60480 & = & 2^6\cdot 3^3\cdot 5\cdot 7 & 34266\cdot 2^6 & \text{} & \text{} & \text{} & \text{} & \text{} \\
 *69120 & = & 2^9\cdot 3^3\cdot 5 & 6616\cdot 2^9 & \text{} & \text{} & \text{} & \text{} & \text{} \\
 \text{} & \text{} & \text{} & \text{} & \text{} & 83160 & = & 2^3\cdot 3^3\cdot 5\cdot 7\cdot 11 & 128 \\
 86400 & = & 2^7\cdot 3^3\cdot 5^2 & 28481\cdot 2^7 & \text{} & \text{} & \text{} & \text{} & \text{} \\
 *103680 & = & 2^8\cdot 3^4\cdot 5 & 23692\cdot 2^8 & \text{} & \text{} & \text{} & \text{} & \text{} \\
 \text{} & \text{} & \text{} & \text{} & \text{} & *110880 & = & 2^5\cdot 3^2\cdot 5\cdot 7\cdot 11 & 144 \\
 *120960 & = & 2^7\cdot 3^3\cdot 5\cdot 7 & 53485\cdot 2^7 & \text{} & \text{} & \text{} & \text{} & \text{} \\
 138240 & = & 2^{10}\cdot 3^3\cdot 5 & 8790\cdot 2^{10} & \text{} & \text{} & \text{} & \text{} & \text{} \\
 161280 & = & 2^9\cdot 3^2\cdot 5\cdot 7 & 17656\cdot 2^9 & \text{} & \text{} & \text{} & \text{} & \text{} \\
 \text{} & \text{} & \text{} & \text{} & \text{} & *166320 & = & 2^4\cdot 3^3\cdot 5\cdot 7\cdot 11 & 160 \\
 *172800 & = & 2^8\cdot 3^3\cdot 5^2 & 42520\cdot 2^8 & \text{} & \text{} & \text{} & \text{} & \text{} \\
 *207360 & = & 2^9\cdot 3^4\cdot 5 & 34026\cdot 2^9 & \text{} & \text{} & \text{} & \text{} & \text{} \\
 \text{} & \text{} & \text{} & \text{} & \text{} & 221760 & = & 2^6\cdot 3^2\cdot 5\cdot 7\cdot 11 & 168 \\
 *241920 & = & 2^8\cdot 3^3\cdot 5\cdot 7 & 80176\cdot 2^8 & \text{} & \text{} & \text{} & \text{} & \text{} \\
 276480 & = & 2^{11}\cdot 3^3\cdot 5 & 11447\cdot 2^{11} & \text{} & \text{} & \text{} & \text{} & \text{} \\
 \text{} & \text{} & \text{} & \text{} & \text{} & *277200 & = & 2^4\cdot 3^2\cdot 5^2\cdot 7\cdot 11 & 180 \\
 311040 & = & 2^8\cdot 3^5\cdot 5 & 103540\cdot 2^8 & \text{} & \text{} & \text{} & \text{} & \text{} \\
 \text{} & \text{} & \text{} & \text{} & \text{} & *332640 & = & 2^5\cdot 3^3\cdot 5\cdot 7\cdot 11 & 192 \\
 *345600 & = & 2^9\cdot 3^3\cdot 5^2 & 61436\cdot 2^9 & \text{} & \text{} & \text{} & \text{} & \text{} \\
 *362880 & = & 2^7\cdot 3^4\cdot 5\cdot 7 & 267219\cdot 2^7 & \text{} & \text{} & \text{} & \text{} & \text{} \\
 *414720 & = & 2^{10}\cdot 3^4\cdot 5 & 47576\cdot 2^{10} & \text{} & \text{} & \text{} & \text{} & \text{} \\
 *483840 & = & 2^9\cdot 3^3\cdot 5\cdot 7 & 116256\cdot 2^9 & \text{} & \text{} & \text{} & \text{} & \text{} \\
 \text{} & \text{} & \text{} & \text{} & \text{} & 498960 & = & 2^4\cdot 3^4\cdot 5\cdot 7\cdot 11 & 200 \\
 552960 & = & 2^{12}\cdot 3^3\cdot 5 & 14652\cdot 2^{12} & \text{} & \text{} & \text{} & \text{} & \text{} \\
 \text{} & \text{} & \text{} & \text{} & \text{} & *554400 & = & 2^5\cdot 3^2\cdot 5^2\cdot 7\cdot 11 & 216 \\
 604800 & = & 2^7\cdot 3^3\cdot 5^2\cdot 7 & 480953\cdot 2^7 & \text{} & \text{} & \text{} & \text{} & \text{} \\
 622080 & = & 2^9\cdot 3^5\cdot 5 & 156278\cdot 2^9 & \text{} & \text{} & \text{} & \text{} & \text{} \\
 \text{} & \text{} & \text{} & \text{} & \text{} & *665280 & = & 2^6\cdot 3^3\cdot 5\cdot 7\cdot 11 & 224 \\
 691200 & = & 2^{10}\cdot 3^3\cdot 5^2 & 86362\cdot 2^{10} & \text{} & \text{} & \text{} & \text{} & \text{} \\
 \text{} & \text{} & \text{} & \text{} & \text{} & *720720 & = & 2^4\cdot 3^2\cdot 5\cdot 7\cdot 11\cdot 13 & 240 \\
 *725760 & = & 2^8\cdot 3^4\cdot 5\cdot 7 & 422932\cdot 2^8 & \text{} & \text{} & \text{} & \text{} & \text{} \\
 829440 & = & 2^{11}\cdot 3^4\cdot 5 & 65018\cdot 2^{11} & \text{} & \text{} & \text{} & \text{} & \text{} \\
 *967680 & = & 2^{10}\cdot 3^3\cdot 5\cdot 7 & 163934\cdot 2^{10} & \text{} & \text{} & \text{} & \text{} & \text{} \\
 \\
 \multicolumn{4}{c}{\emph{Recursively super-abundant but}} 	&&		\\
 \multicolumn{4}{c}{\emph{not recursively highly composite}} 	&&  		\\
  *181440 & = & 2^6\cdot 3^4\cdot 5\cdot7 				&&  		
\end{array}
$
\end{center}
\end{small}
\begin{footnotesize}

\end{footnotesize}
\end{document}